\documentclass[a4paper,notitlepage]{article}%
\usepackage{amssymb}
\usepackage{graphicx}
\usepackage{amsmath}
\usepackage{amsthm}
\usepackage{amsfonts}%
\providecommand{\U}[1]{\protect\rule{.1in}{.1in}}

\providecommand{\U}[1]{\protect\rule{.1in}{.1in}}

\theoremstyle{plain}
\newtheorem{theorem}{Theorem}[section]

\newtheorem{corollary}[theorem]{Corollary}

\newtheorem{lemma}[theorem]{Lemma}

\theoremstyle{definition}

\newtheorem{remark}[theorem]{Remark}
\numberwithin{equation}{section}
\numberwithin{theorem}{section}

\begin{document}

\title{Positive solutions for concave-convex type problems for the one-dimensional
$\phi$-Laplacian\thanks{2020 \textit{Mathematics Subject Clasification}.
34B15; 34B18.} 
\thanks{Partially supported by Secyt-UNC 33620180100016CB.} }
\author{Uriel Kaufmann, Leandro Milne\thanks{\textit{E-mail addresses. }
kaufmann@mate.uncor.edu (Uriel Kaufmann),
leandro.milne@unc.edu.ar (Leandro Milne, Corresponding Author).}
\and \noindent\\{\small FAMAF, Universidad Nacional de C\'{o}rdoba, (5000) C\'{o}rdoba,
Argentina}}
\maketitle

\begin{abstract}
Let $\Omega=(a,b)\subset\mathbb{R}$, $0\leq m,n\in L^{1}(\Omega)$,
$\lambda,\mu>0$ be real parameters, and $\phi:\mathbb{R}\rightarrow\mathbb{R}$
be an odd increasing homeomorphism. In this paper we consider the existence of
positive solutions for problems of the form
\[%
\begin{cases}
-\phi\left(  u^{\prime}\right)  ^{\prime}=\lambda m(x)f(u)+\mu n(x)g(u) &
\text{ in }\Omega,\\
u=0 & \text{ on }\partial\Omega,
\end{cases}
\]
where $f,g:[0,\infty)\rightarrow\lbrack0,\infty)$ are continuous functions
which are, roughly speaking, sublinear and superlinear with respect to $\phi$,
respectively. Our assumptions on $\phi$, $m$ and $n$ are substantially weaker
than the ones imposed in previous works. The approach used here combines the
Guo-Krasnoselski\u{\i}\ fixed-point theorem and the sub-supersolutions method
with some estimates on related nonlinear problems.

\end{abstract}

\section{Introduction}

Let $\Omega=\left(  a,b\right)  \subset\mathbb{R}$, $m,n\in L^{1}\left(
\Omega\right)  $ and $\lambda,\mu>0$ be a real parameters. In this article we
consider problems of the form
\begin{equation}
\left\{
\begin{array}
[c]{ll}%
-\phi\left(  u^{\prime}\right)  ^{\prime}=\lambda m\left(  x\right)  f\left(
u\right)  +\mu n(x)g(u) & \text{in }\Omega,\\
u=0 & \text{on }\partial\Omega,
\end{array}
\right.  \label{problema concavo convexo}%
\end{equation}
where $\phi:\mathbb{R\rightarrow R}$ is an odd increasing homeomorphism and
$f,g:\left[  0,\infty\right)  \rightarrow\left[  0,\infty\right)  $ are
continuous functions which are, roughly speaking, sublinear and superlinear
with respect to $\phi$, respectively.
When the nonlinearities $f$ and $g$ are concave and convex, the problem
\eqref{problema concavo convexo} with $\phi(x)=x$ was first studied by
Ambrosetti, Brezis and Cerami in their celebrated paper
\cite{ambrosettibreziscerami}. More precisely, in that article the authors
studied the $N$-dimensional problem
\begin{equation}
\left\{
\begin{array}
[c]{ll}%
-\Delta u=\lambda u^{q}+u^{p} & x\in\Omega,\\
u>0 & x\in\Omega,\\
u=0 & x\in\partial\Omega,
\end{array}
\right.  \label{2}%
\end{equation}
with $0<q<1<p$ and $\Omega$ a bounded domain in $\mathbb{R}^{N}$. They proved
the following facts: there exists $\Lambda>0$ such that: if $\lambda
\in(0,\Lambda)$ then \eqref{2} has at least two positive solutions, if
$\lambda=\Lambda$ there is at least one positive solution, and if
$\lambda>\Lambda$ then there are no positive solutions.

Several authors have studied generalizations of \eqref{2}, see for instance
\cite{ambrosettigarciaperal1996, papageorgiou_smyrlis2016, sanchezubilla2000}
and their references, where the corresponding problem for the $p$-Laplacian is
considered. Also, in \cite{papageorgiou_winkert2016} the authors have treated
the $N$-dimensional problem for the $\phi$-Laplacian operator.

Regarding the one-dimensional $\phi$-Laplacian problem that we will deal with
in this article, Wang in \cite[Theorem 1.2]{wang} and \cite[Theorem
1.2]{wang2} studied the case
$m=n\geq0$, $m\not \equiv 0$ on any subinterval in $\Omega$, $m\in
C(\overline{\Omega})$ and $\lambda=\mu$. In these papers it is proved that
there exist $\lambda_{0},\lambda_{1}>0$ such that if $\lambda\in(0,\lambda
_{0})$, then \eqref{problema concavo convexo} has at least two positive
solutions; and if $\lambda>\lambda_{1}$, then there are no positive solutions.
Let us note that the hypothesis on $\phi$ imposed in \cite{wang, wang2} are
much stronger than the ones that we shall require here. More precisely, Wang assumes

\begin{enumerate}
\item[$\left(  \Phi\right)  $] There exist increasing homeomorphisms $\psi
_{1},\psi_{2}:\left[  0,\infty\right)  \rightarrow\left[  0,\infty\right)  $
such that $\psi_{1}\left(  t\right)  \phi\left(  x\right)  \leq\phi\left(
tx\right)  \leq\psi_{2}\left(  t\right)  \phi\left(  x\right)  $ for all
$t,x>0.$
\end{enumerate}

\noindent On other hand, \eqref{problema concavo convexo} is also considered
in \cite{lee2020} with $m=n\geq0$, $m\not \equiv 0$ on any subinterval in
$\Omega$ and $\lambda=\mu$ like in \cite{wang,wang2}. However, the regularity
assumptions for $m $ allow some $m \in L_{loc}^{1}(\Omega) $. Regarding the
hypothesis on $\phi$ they require that

\begin{enumerate}
\item[$\left(  \Phi^{\prime}\right)  $] There exist an increasing
homeomorphism $\psi_{1}:\left[  0,\infty\right)  \rightarrow\left[
0,\infty\right)  $ and a function $\psi_{2}:\left[  0,\infty\right)
\rightarrow\left[  0,\infty\right)  $ such that $\psi_{1}\left(  t\right)
\phi\left(  x\right)  \leq\phi\left(  tx\right)  \leq\psi_{2}\left(  t\right)
\phi\left(  x\right)  $ for all $t,x>0.$
\end{enumerate}

\noindent The authors prove that there exist $\lambda_{1}\geq\lambda_{0}>0$
such that \eqref{problema concavo convexo} has at least two positive solutions
for $\lambda\in(0,\lambda_{0})$, one positive solution for $\lambda\in
\lbrack\lambda_{0},\lambda_{1}]$, and no positive solution for $\lambda
>\lambda_{1}$.

In this article, employing the method of sub and supersolutions and the
Guo-Krasnoselski\u{\i} fixed-point theorem along with some estimates for
related problems, we shall prove that there are at least two positive
solutions for $\lambda\approx0$, under much weaker assumptions on $\phi$, $m$
and $n$. Moreover, as a consequence of Theorem
\ref{equivalencias de phi y phi'} we shall see that $(\Phi)$ and
$(\Phi^{\prime})$ are in fact equivalent.

To be more precise, let us introduce the following hypothesis.

\begin{enumerate}
\item[(F)] There exist $c_{0},t_{0},q>0$ such that
\begin{equation}
f(t)\geq c_{0}t^{q}\text{ for all }t\in\lbrack0,t_{0}]\quad\text{and}\quad
\lim_{t\rightarrow0^{+}}\frac{t^{q}}{\phi(t)}=\infty. \label{(F)}%
\end{equation}

\item[(G1)] There exist $c_{1},t_{1},r_{1}>0$ such that
\begin{equation}
g(t)\leq c_{1}t^{r_{1}}\text{ for all }t\in\lbrack0,t_{1}]\quad\text{and}%
\quad\lim_{t\rightarrow0^{+}}\frac{t^{r_{1}}}{\phi(t)}=0. \label{G1}%
\end{equation}

\item[(G2)] There exist $c_{2},t_{2},r_{2}>0$ such that
\begin{equation}
g(t)\geq c_{2}t^{r_{2}}\text{ for all }t\geq t_{2}\quad\text{and}\quad
\lim_{t\rightarrow\infty}\frac{t^{r_{2}}}{\phi(t)}=\infty. \label{(G2)}%
\end{equation}

\end{enumerate}

\noindent Note that when $\phi(t)=\left\vert t\right\vert ^{p-2}t$,
$f(u)=u^{q}$ and $g(u)=u^{r}$, the limits in (F) and (G1) are satisfied if and
only if $0<q<p-1<r$. Let us set $\mathcal{C}_{0}^{1}(\overline{\Omega
}):=\{u\in\mathcal{C}^{1}(\overline{\Omega}):u=0$ on $\partial\Omega\}$ and%
\[
\mathcal{P}^{\circ}:=\left\{  u\in\mathcal{C}_{0}^{1}(\overline{\Omega
}):u>0\text{ in }\Omega\text{\ and }u^{\prime}\left(  b\right)  <0<u^{\prime
}\left(  a\right)  \right\}  .
\]

Our main result is the following theorem:\smallskip

\begin{theorem}
\strut\label{teorema principal prob conc convex}Let $0\leq m,n\in L^{1}\left(
\Omega\right)  $.

\begin{enumerate}
\item[(I)] \label{item de existencia con sub y supersols prob conc convex}
Assume that $m\not \equiv 0$ and (F) and (G1) hold. Then for all $\mu>0$ there
exists $\lambda_{0}(\mu)>0$ such that (\ref{problema concavo convexo}) has a
solution $u_{\lambda}\in\mathcal{P}^{\circ}$ for all $0<\lambda<\lambda
_{0}(\mu)$. Moreover, the solutions $u_{\lambda}$ can be chosen such that
\begin{equation}
\lim_{\lambda\rightarrow0^{+}}\left\Vert u_{\lambda}\right\Vert _{\infty}=0.
\label{algunas soluciones tienden a cero cuando lambda tiende a cero}%
\end{equation}

\item[(II)] \label{item de existencia usando punto fijo prob conc convex}
Assume that $n\not \equiv 0$ and (G1) and (G2) hold. Then for all $\mu>0$
there exists $\lambda_{1}(\mu)>0$ such that (\ref{problema concavo convexo})
has a solution $v_{\lambda}\in\mathcal{P}^{\circ}$ for all $0<\lambda
<\lambda_{1}(\mu)$. Furthermore, there exists $\rho>0$ such that $\left\Vert
v_{\lambda}\right\Vert _{\infty}>\rho$ for all $0<\lambda<\lambda_{1}(\mu).$

\item[(III)] \label{item de conjunto de soluciones abierto} Assume that
$\{\lambda>0:\eqref{problema concavo convexo}\ \text{has a solution
in}\ \mathcal{P}^{\circ}\}\neq\emptyset$ and (F) holds for all $t_{0}>0$. Let
\[
\Lambda:=\sup\{\lambda>0:\eqref{problema concavo convexo}\ \text{has a
solution in}\ \mathcal{P}^{\circ}\}.
\]
Then, for $0<\lambda<\Lambda$ $\eqref{problema concavo convexo}$ has at least
one solution in $\mathcal{P}^{\circ}$.
\end{enumerate}
\end{theorem}

As an immediate consequence of the above theorem we have the
following\smallskip

\begin{corollary}
Let $\mu>0$ and $0\leq m,n\in L^{1}(\Omega)$ with $m,n\not \equiv 0$. Assume
that (F), (G1) and (G2) hold. Then \eqref{problema concavo convexo} has at
least two solutions in $\mathcal{P}^{\circ}$ for $\lambda\approx0$.
\end{corollary}

The rest of the paper is organized as follows. In the next section we state
some necessary facts about nonlinear problems involving the $\phi$-Laplacian,
and in Section \ref{main results} we prove our main results. Finally, in
Section \ref{section comments} we introduce some concepts about Orlicz spaces
indices which we use to prove Theorem \ref{equivalencias de phi y phi'} (and,
in particular, the equivalence of $(\Phi)$ and $(\Phi^{\prime})$), and at the
end of the section we give several examples of functions $\phi$ illustrating
our conditions and their relations with the ones used in the previous works.
Let us mention that all the $\phi$'s constructed in Example (e) satisfy
conditions (F), (G1) and (G2) but do not fulfill condition $\left(
\Phi\right)  $.

\section{Preliminaries}

Let $\phi:\mathbb{R\rightarrow R}$ be an odd increasing homeomorphism. We
start considering problems of the form
\begin{equation}
\left\{
\begin{array}
[c]{ll}%
-\phi\left(  v^{\prime}\right)  ^{\prime}=h\left(  x\right)  & \text{in
}\Omega,\\
v=0 & \text{on }\partial\Omega.
\end{array}
\right.  \label{g}%
\end{equation}
It is well known that for all $h\in L^{1}(\Omega)$, (\ref{g}) possesses a
unique solution $v\in\mathcal{C}_{0}^{1}(\overline{\Omega})$ such that
$\phi\left(  v^{\prime}\right)  $ is absolutely continuous and that the
equation holds pointwise $a.e.$ $x\in\Omega$. Furthermore, the solution
operator $\mathcal{S}_{\phi}\colon L^{1}(\Omega)\rightarrow\mathcal{C}%
^{1}(\overline{\Omega})$ is continuous and nondecreasing, see \cite[Lemma
2.1]{dang} and \cite[Lemma 2.2]{JMAA2018}.

We need now to introduce some notation. For $0\leq h\in L^{1}(\Omega)$ with
$h\not \equiv 0$, set
\begin{align*}
\mathcal{A}_{h}  &  :=\left\{  x\in\Omega:h\left(  y\right)  =0\text{
}a.e.\text{ }y\in\left(  a,x\right)  \right\}  ,\\
\mathcal{B}_{h}  &  :=\left\{  x\in\Omega:h\left(  y\right)  =0\text{
}a.e.\text{ }y\in\left(  x,b\right)  \right\}  ,
\end{align*}
and
\begin{gather}
\alpha_{h}:=\left\{
\begin{array}
[c]{ll}%
\sup\mathcal{A}_{h} & \text{if }\mathcal{A}_{h}\not =\emptyset,\\
a & \text{if }\mathcal{A}_{h}=\emptyset,
\end{array}
\right.  \quad\beta_{h}:=\left\{
\begin{array}
[c]{ll}%
\inf\mathcal{B}_{h} & \text{if }\mathcal{B}_{h}\not =\emptyset,\\
b & \text{if }\mathcal{B}_{h}=\emptyset,
\end{array}
\right. \nonumber\\
\underline{\theta}_{h}:=\min\left\{  \frac{1}{\beta_{h}-a},\frac{1}%
{b-\alpha_{h}}\right\}  ,\quad\overline{\theta}_{h}:=\frac{\alpha_{h}%
+\beta_{h}}{2}. \label{tita}%
\end{gather}
We observe that $\underline{\theta}_{h}$ is well defined because
$h\not \equiv 0$, and $\alpha_{h}<\beta_{h}$ (and so, $\overline{\theta}%
_{h}\in\left(  \alpha_{h},\beta_{h}\right)  $). We also write
\[
\delta_{\Omega}\left(  x\right)  :=\text{dist}\left(  x,\partial\Omega\right)
=\min\left(  x-a,b-x\right)  \text{.}%
\]

We shall utilize the following estimates on several occasions in the sequel.
For the proof, see \cite[Lemma 2.3 and (2.6)]{JMAA2018} and \cite[Corollary
2.2]{JMAA2019}.

\begin{lemma}
\label{lemma cotas}Let $0\leq h\in L^{1}(\Omega)$ with $h\not \equiv 0$.

\begin{enumerate}
\item[(i)] In $\overline{\Omega}$ it holds that
\begin{gather}
\underline{\theta}_{h}\min\left\{  \int_{a}^{\overline{\theta}_{h}}\phi
^{-1}\left(  \int_{y}^{\overline{\theta}_{h}}h\right)  dy,\int_{\overline
{\theta}_{h}}^{b}\phi^{-1}\left(  \int_{\overline{\theta}_{h}}^{y}h\right)
dy\right\}  \delta_{\Omega}\nonumber\\
\leq\mathcal{S}_{\phi}\left(  h\right)  \leq\phi^{-1}\left(  \int_{a}%
^{b}h\right)  \delta_{\Omega}. \label{des1}%
\end{gather}

\item[(ii)] In $\overline{\Omega}$ it holds that
\begin{equation}
\mathcal{S_{\phi}}(h)\geq\underline{\theta}_{h}\left\Vert \mathcal{S}_{\phi
}(h)\right\Vert _{\infty}\delta_{\Omega}.
\label{desigualdad para probar inciso iii del concavo convexo}%
\end{equation}

\item[(iii)] For $M>0$ there exists $c>0$ not depending on $M$ such that in
$\overline{\Omega}$ it holds that%
\begin{equation}
\min\left\{  \int_{a}^{\overline{\theta}_{h}}\phi^{-1}\left(  \int
_{y}^{\overline{\theta}_{h}}Mh\right)  dy,\int_{\overline{\theta}_{h}}^{b}%
\phi^{-1}\left(  \int_{\overline{\theta}_{h}}^{y}Mh\right)  dy\right\}  \geq
c\phi^{-1}(cM). \label{des3}%
\end{equation}

\end{enumerate}
\end{lemma}

\noindent Observe that, since $\overline{\theta}_{h}\in\left(  \alpha
_{h},\beta_{h}\right)  $, the constant that appears in the first term of the
inequalities in (\ref{des1}) is strictly positive. Note also that, since
$\underline{\theta}_{h}\left\Vert \delta_{\Omega}\right\Vert _{\infty}\geq
1/2$, using the lower bound of (\ref{des1}) and taking into account the
monotonicity of the infinite norm we get
\begin{equation}
\frac{1}{2}\min\left\{  \int_{a}^{\overline{\theta}_{h}}\phi^{-1}\left(
\int_{y}^{\overline{\theta}_{h}}h\right)  dy,\int_{\overline{\theta}_{h}}%
^{b}\phi^{-1}\left(  \int_{\overline{\theta}_{h}}^{y}h\right)  dy\right\}
\leq\left\Vert \mathcal{S}_{\phi}\left(  h\right)  \right\Vert _{\infty}.
\label{des2}%
\end{equation}
Observe also that for $h$ as in Lemma \ref{lemma cotas} $\mathcal{S}_{\phi
}\left(  h\right)  \in\mathcal{P}^{\circ}$.

Let $h:\Omega\times\mathbb{R}\rightarrow\mathbb{R}$ be a Carath\'{e}odory
function (that is, $h\left(  x,\cdot\right)  $ is continuous for $a.e.$
$x\in\Omega$ and $h\left(  \cdot,\xi\right)  $ is measurable for all $\xi
\in\mathbb{R}$). We now consider problems of the form
\begin{equation}
\left\{
\begin{array}
[c]{ll}%
-\phi\left(  u^{\prime}\right)  ^{\prime}=h\left(  x,u\right)  & \text{in
}\Omega,\\
u=0 & \text{on }\partial\Omega.
\end{array}
\right.  \label{noli}%
\end{equation}
We shall say that $v\in\mathcal{C}(\overline{\Omega})$ is a
\textit{subsolution }of (\ref{noli}) if there exists a finite set
$\Sigma\subset\Omega$ such that $\phi(v^{\prime})\in AC_{loc}(\overline
{\Omega}\,\backslash\,\Sigma),$ $v^{\prime}(\tau^{+}):=\lim_{x\rightarrow
\tau^{+}}v^{\prime}(x)\in\mathbb{R}$, $v^{\prime}(\tau^{-}):=\lim
_{x\rightarrow\tau^{-}}v^{\prime}(x)\in\mathbb{R}$ for each $\tau\in\Sigma,$
and
\begin{equation}
\left\{
\begin{array}
[c]{ll}%
-\phi\left(  v^{\prime}\right)  ^{\prime}\leq h\left(  x,v\left(  x\right)
\right)  & a.e.\text{ }x\in\Omega,\\
v\leq0\text{ on }\partial\Omega, & v^{\prime}(\tau^{-})<v^{\prime}(\tau
^{+})\text{ for each }\tau\in\Sigma.
\end{array}
\right.  \label{subi}%
\end{equation}
If the inequalities in (\ref{subi}) are inverted, we shall say that $v$ is a
\textit{supersolution} of (\ref{noli}).\medskip

For the sake of completeness, we state an existence result in the presence of
well-ordered sub and supersolutions, and a particular case of the well-known
Guo-Krasnoselski\u{\i} fixed-point theorem (for a proof, see e.g.
\cite[Theorem 7.16]{rach} and \cite[Theorem 2.3.4]{guo}, respectively).

\begin{lemma}
\label{subsup}Let $v$ and $w$ be sub and supersolutions respectively of
(\ref{noli}) such that $v\leq w$ in $\Omega$. Suppose there exists $g\in
L^{1}\left(  \Omega\right)  $ such that
\[
\left\vert h\left(  x,\xi\right)  \right\vert \leq g\left(  x\right)
\quad\text{for }a.e.\text{ }x\in\Omega\text{ and all }\xi\in\left[  v\left(
x\right)  ,w\left(  x\right)  \right]  .
\]
Then there exists $u\in\mathcal{C}_{0}^{1}(\overline{\Omega})$ solution of
(\ref{noli}) with $v\leq u\leq w$ in $\Omega$.
\end{lemma}

\begin{lemma}
\label{teorema de punto fijo de conos en expansion y compresion} Let $X $ be a
Banach space and let $K $ be a cone in X. Let $\Omega_{1} , \Omega_{2} \subset
X $ be two open sets with $0 \in\Omega_{1} $ and $\Omega_{1}\subset\Omega_{2}
$. Suppose that $T : K \cap(\Omega_{2} \setminus\Omega_{1} ) \rightarrow K $
is a completely continuous operator and
\begin{align}
\left\Vert Tv \right\Vert \ge\left\Vert v \right\Vert , \quad\text{for} \ v
\in K\cap\partial\Omega_{2},\nonumber\\
\left\Vert Tv \right\Vert \le\left\Vert v \right\Vert , \quad\text{for} \ v
\in K\cap\partial\Omega_{1}.\nonumber
\end{align}
Then, $T$ has a fixed point in $K \cap(\Omega_{2} \setminus\Omega_{1} ) $.
\end{lemma}

\section{Proof of the main results}

\label{main results}

\subsection{Proof of item (I)}

We start this section with two lemmas concerning sub and supersolutions that
shall be used to prove item (I) of Theorem
\ref{teorema principal prob conc convex}.

\begin{lemma}
\label{lema supersol del problema concavo convexo} Let $m,n\in L^{1}(\Omega)$
such that $0\not \equiv m+n\geq0$. Assume that (G1) holds. Then for all
$\mu>0$ there exists $\lambda_{0}(\mu)>0$ such that for each $0<\lambda
<\lambda_{0}(\mu)$ there exists $w_{\lambda}\in\mathcal{P}^{\circ}$
supersolution of (\ref{problema concavo convexo}). Moreover,
\begin{equation}
\lim_{\lambda\rightarrow0^{+}}\left\Vert w_{\lambda}\right\Vert _{\infty}=0.
\label{supsol comverge to zero}%
\end{equation}

\end{lemma}

\textit{Proof}. Let $c_{1},t_{1},r_{1}$ be given by (G1). Let us define
$c_{\Omega}:=\max_{\overline{\Omega}}\delta_{\Omega}$. By the continuity of
$\phi^{-1}$ and the fact that $\phi^{-1}(0)=0$, there exists $K_{0}>0$ such
that
\begin{equation}
\phi^{-1}(\kappa\int_{a}^{b}m(s)+n(s)ds)\leq\frac{t_{1}}{c_{\Omega}}%
\quad\text{for all }\kappa\leq K_{0}%
.\label{lema supersol del problema concavo convexo 0}%
\end{equation}
We observe that by the second condition on (\ref{G1}), for $\rho>0$ fixed we
have
\begin{equation}
\lim_{t\rightarrow0^{+}}\frac{[\phi^{-1}(\rho t)]^{r_{1}}}{t}%
=0.\label{lema supersol del problema concavo convexo 1}%
\end{equation}

We now define
\[
\epsilon:=\frac{1}{c_{1}\mu c_{\Omega}^{r_{1}}},\quad\rho:=\int_{a}%
^{b}m(s)+n(s)ds.
\]
We can deduce from (\ref{lema supersol del problema concavo convexo 1}) that
there exists $K_{1}=K_{1}(\epsilon,\rho)>0$ such that
\begin{equation}
\lbrack\phi^{-1}(\kappa\rho)]^{r_{1}}\leq\kappa\epsilon\text{ for all }%
\kappa\leq K_{1}.\label{lema supersol del problema concavo convexo 3}%
\end{equation}
Let $C=\max_{[0,t_{1}]}f(t)$ and choose $\lambda_{0}>0$ such that
\begin{equation}
\lambda_{0}C\leq\min\{K_{0},K_{1}%
\}.\label{lema supersol del problema concavo convexo 5}%
\end{equation}
Also, for each $0<\lambda<\lambda_{0}$, pick $\kappa_{\lambda}$ such that
\begin{equation}
\lambda C\leq\kappa_{\lambda}\leq\min\{K_{0},K_{1}%
\},\label{lema supersol del problema concavo convexo 4}%
\end{equation}
and for such $\kappa_{\lambda}$ define $w_{\lambda}:=\mathcal{S}_{\phi}%
(\kappa_{\lambda}(m+n)).$ Since $\kappa_{\lambda}\leq K_{0}$, the upper bound
in \eqref{des1} and (\ref{lema supersol del problema concavo convexo 0}) tell
us that $\left\Vert w_{\lambda}\right\Vert _{\infty}\leq t_{1}.$ Taking into
account (\ref{lema supersol del problema concavo convexo 3}),
(\ref{lema supersol del problema concavo convexo 5}) and
(\ref{lema supersol del problema concavo convexo 4}), employing (G1) and the
upper bound in \eqref{des1} we deduce that
\begin{gather*}
\lambda m(x)f(w_{\lambda})+\mu n(x)g(w_{\lambda})\\
\leq\lambda m(x)C+c_{1}\mu n(x)w_{\lambda}^{r_{1}}\\
\leq\kappa_{\lambda}m(x)+c_{1}\mu n(x)\left[  \phi^{-1}(\kappa_{\lambda}%
\int_{a}^{b}m(s)+n(s)ds)\delta_{\Omega}\right]  ^{r_{1}}\\
\leq\kappa_{\lambda}(m(x)+n(x))=-\phi(w_{\lambda}^{\prime})^{\prime}%
\quad\text{in }\Omega,
\end{gather*}
and hence $w_{\lambda}$ is a supersolution of (\ref{problema concavo convexo}).

In order to prove \eqref{supsol comverge to zero}, we choose $\kappa_{\lambda
}$ satisfying \eqref{lema supersol del problema concavo convexo 4} and such
that $\kappa_{\lambda}\rightarrow0$ when $\lambda\rightarrow0^{+}$. Hence,
using the second inequality \eqref{des1} we get that
\[
0\leq w_{\lambda}(x)=\mathcal{S}_{\phi}(\kappa_{\lambda}(m+n))\leq\phi
^{-1}\left(  \int_{a}^{b}\kappa_{\lambda}(m+n)\right)  \delta_{\Omega
}(x)\rightarrow0
\]
uniformly in $\overline{\Omega}$ when $\lambda\rightarrow0^{+}$. Thus,
$\lim_{\lambda\rightarrow0^{+}}\Vert w_{\lambda}\Vert_{\infty}=0$. \qed

\begin{lemma}
\label{lema subsol del problema concavo convexo} Let $0\leq m,n\in
L^{1}(\Omega)$ with $m\not \equiv 0.$ Assume that (F) holds. Then for all
$\lambda,\mu>0$ (\ref{problema concavo convexo}) has a subsolution
$v\in\mathcal{P}^{\circ}$.
\end{lemma}

\textit{Proof.}
Let $\lambda,\mu>0$ and let $c_{0},t_{0},q$ be given by (F). Recall that
$c_{\Omega}:=\max_{\overline{\Omega}}\delta_{\Omega}.$ Since $\phi^{-1}$ is
continuous and $\phi^{-1}(0)=0$, there exists $\varepsilon_{0}>0$ such that
\begin{equation}
\phi^{-1}(\varepsilon\int_{a}^{b}m(s)\delta_{\Omega}^{q}(s)ds)\leq\frac{t_{0}%
}{c_{\Omega}}\quad\text{for all }\varepsilon\leq\varepsilon_{0}%
.\label{lema subsol del problema concavo convexo 0}%
\end{equation}
By the second condition in (\ref{(F)}), for $\rho>0$ fixed
\begin{equation}
\lim_{t\rightarrow0^{+}}\frac{[\phi^{-1}(\rho t)]^{q}}{t}=\infty
.\label{lema subsol del problema concavo convexo 1}%
\end{equation}
Let us define
\[
M:=\frac{1}{\lambda c_{0}c^{q}},
\]
where $c$ is the constant in \eqref{des3} with $h=m\delta_{\Omega}^{q}$. It
follows from (\ref{lema subsol del problema concavo convexo 1}) that there
exists $\varepsilon_{1}=\varepsilon_{1}(M,\rho)$ such that
\begin{equation}
\lbrack\phi^{-1}(\varepsilon\rho)]^{q}\geq M\varepsilon\text{ for all
}\varepsilon\leq\varepsilon_{1}%
.\label{lema subsol del problema concavo convexo 2}%
\end{equation}

Let us choose
\begin{equation}
0<\varepsilon<\min\{\varepsilon_{0},\varepsilon_{1}%
\}\label{lema subsol del problema concavo convexo 3}%
\end{equation}
and for such $\varepsilon$ define $v:=\mathcal{S}_{\phi}(\varepsilon
m\delta_{\Omega}^{q}).$ Since $\varepsilon\leq\varepsilon_{0}$, the upper
bound of Lemma \ref{lemma cotas} and
(\ref{lema subsol del problema concavo convexo 0}) tell us that $\left\Vert
v\right\Vert _{\infty}\leq t_{0}$. Consequently, taking into account
(\ref{lema subsol del problema concavo convexo 2}) and
(\ref{lema subsol del problema concavo convexo 3}), employing (F) and
(\ref{des3}) we deduce that
\begin{gather*}
\lambda m(x)f(v)+\mu n(x)g(v)\geq\lambda c_{0}m(x)v^{q}\\
\geq\lambda c_{0}m(x)[c\phi^{-1}(c\varepsilon)\delta_{\Omega}]^{q}%
\geq\varepsilon m(x)\delta_{\Omega}^{q}\quad\text{in }\Omega.
\end{gather*}
In other words, $v$ is a subsolution of (\ref{problema concavo convexo}). \qed
\medskip

\textit{Proof of Theorem 1.1 (I)}.
Given $\mu>0$, let $\lambda_{0}(\mu)$ be as in Lemma
\ref{lema supersol del problema concavo convexo}. For $0<\lambda<\lambda
_{0}(\mu)$, let $w_{\lambda}\in\mathcal{P}^{\circ}$ be a supersolution
provided by the aforementioned lemma, and let $v_{\lambda}\in\mathcal{P}%
^{\circ}$ be a subsolution given by Lemma
\ref{lema subsol del problema concavo convexo} with $\varepsilon_{\lambda}$
chosen such that $\varepsilon_{\lambda}m(x)\delta_{\Omega}^{q}(x)\leq
\kappa_{\lambda}(m(x)+n(x))$ for $a.e.x\in\Omega$. It follows that
$v_{\lambda},w_{\lambda}$ are a pair of well-ordered sub and supersolutions of
\eqref{problema concavo convexo}. Hence, Lemma \ref{subsup} gives a solution
of \eqref{problema concavo convexo} $u_{\lambda}\in\mathcal{P}^{\circ}$.
Moreover,
\eqref{algunas soluciones tienden a cero cuando lambda tiende a cero} follows
from \eqref{supsol comverge to zero}. \qed

\subsection{Proof of item (II)}

\textit{Proof of Theorem 1.1 (II)}. We shall use Lemma
\ref{teorema de punto fijo de conos en expansion y compresion} with the
operator
\[
Tv:=\mathcal{S}_{\phi}(\lambda m(x)f(v)+\mu n(x)g(v)),
\]
the cone
\[
\mathcal{K}:=\{v\in C(\overline{\Omega}):v\geq\underline{\theta}_{n}\left\Vert
v\right\Vert _{\infty}\delta_{\Omega}\}
\]
($\underline{\theta}_{n}$ as in (\ref{tita})) and the open balls $B_{R}%
(0),B_{\rho}(0)\subset C(\overline{\Omega})$ with $0<\rho<R.$ Observe that
$C_{0}^{1}(\overline{\Omega})\cap(\mathcal{K}\setminus\{0\})\subset
\mathcal{P}^{\circ}$ and any fixed point of $T$ belong to $C^{1}(\overline{\Omega})$.  

Let $c_{2},t_{2}$ and $r_{2}$ be given by (G2). We consider the function
$h:=c_{2}\mu\left(  \underline{\theta}_{n}\right)  ^{r_{2}}n\delta_{\Omega
}^{r_{2}}$. Taking into account \eqref{des3}, we can find $c=c(\mu)>0$ such
that for all $M>0$
\begin{equation}
\min\left\{  \int_{a}^{\overline{\theta}_{n}}\phi^{-1}(M\int_{y}%
^{\overline{\theta}_{n}}h)dy,\int_{\overline{\theta}_{n}}^{b}\phi^{-1}%
(M\int_{\overline{\theta}_{n}}^{y}h)dy\right\}  \geq c\phi^{-1}(cM).
\label{eq1 prueba de existencia con pto fijo conc convex}%
\end{equation}
On other hand, the second condition in (G2) is equivalent to
\[
\lim_{t\rightarrow\infty}\frac{\phi^{-1}(\rho t^{r_{2}})}{t}=\infty
\]
for all fixed $\rho>0$, and then there exists $\overline{t}>0$ such that
\begin{equation}
\phi^{-1}(ct^{q_{2}})\geq\frac{2t}{c}\ \text{for all}\ t\geq\overline{t}.
\label{desigualdad usando G2}%
\end{equation}

Let us fix $R>\max\{t_{2},\overline{t}\}$. Taking into account that
$\mathcal{S}_{\phi}$ and $\phi^{-1}$ are nondecreasing, the inequality
\eqref{des2}, (G2), \eqref{eq1 prueba de existencia con pto fijo conc convex}
and \eqref{desigualdad usando G2} we obtain that for $v\in\mathcal{K}%
\cap\partial B_{R}(0)$,
\begin{align*}
\left\Vert Tv\right\Vert _{\infty}=  &  \left\Vert \mathcal{S}_{\phi}(\lambda
m(x)f(v))+\mu n(x)g(v))\right\Vert \geq\left\Vert S_{\phi}(\mu
n(x)g(v))\right\Vert _{\infty}\\
\geq &  \frac{1}{2}\min\left\{  \int_{a}^{\overline{\theta}_{n}}\phi
^{-1}\left(  \int_{y}^{\overline{\theta}_{n}}\mu ng(v)\right)  dy,\int
_{\overline{\theta}_{n}}^{b}\phi^{-1}\left(  \int_{\overline{\theta}_{n}}%
^{y}\mu ng(v)\right)  dy\right\} \\
\geq &  \frac{1}{2}\min\left\{  \int_{a}^{\overline{\theta}_{n}}\phi
^{-1}\left(  c_{2}\mu\int_{y}^{\overline{\theta}_{n}}nv^{r_{2}}\right)
dy,\int_{\overline{\theta}_{n}}^{b}\phi^{-1}\left(  c_{2}\mu\int
_{\overline{\theta}_{n}}^{y}nv^{r_{2}}\right)  dy\right\} \\
\geq &  \frac{1}{2}\min\left\{  \int_{a}^{\overline{\theta}_{n}}\phi
^{-1}\left(  c_{2}\mu\left(  {\underline{\theta}_{n}}\left\Vert v\right\Vert
_{\infty}\right)  ^{r_{2}}\int_{y}^{\overline{\theta}_{n}}n\delta_{\Omega
}^{r_{2}}\right)  dy,\right. \\
&  \left.  \int_{\overline{\theta}_{n}}^{b}\phi^{-1}\left(  c_{2}\mu\left(
{\underline{\theta}_{n}}\left\Vert v\right\Vert _{\infty}\right)  ^{r_{2}}%
\int_{\overline{\theta}_{n}}^{y}n\delta_{\Omega}^{r_{2}}\right)  dy\right\} \\
\geq &  \frac{1}{2}{c}\phi^{-1}({c}\left\Vert v\right\Vert _{\infty}^{r_{2}%
})\\
\geq &  \left\Vert v\right\Vert _{\infty}.
\end{align*}
That is, $\left\Vert Tv\right\Vert _{\infty}\geq\left\Vert v\right\Vert
_{\infty}$ for such $v.$

On other side, let $N:=c_{1}\int_{a}^{b}n.$ The second condition in (G1)
implies that there exists $\underline{t}>0$ such that $\phi(t/c_{\Omega})>\mu
Nt^{r_{1}}$ for all $t\in(0,\underline{t})$. Set $C:=\max_{[0,R]}f(t)$ and
$M:=\int_{a}^{b}m$. Let $0<\rho<\min\{\underline{t},R/2,t_{1}\}$ be fixed and
define
\begin{equation}
\lambda_{1}:=\frac{\phi(\rho/c_{\Omega})-\mu N\rho^{r_{1}}}{MC}.
\label{conc convex definicion de lambda1}%
\end{equation}
Note that $\lambda_{1}>0$ by our election of $\underline{t}$.

Now, taking into account (\ref{des1}), (G1),
(\ref{conc convex definicion de lambda1}) and the monotonicity of $\phi^{-1}$
we see for $0<\lambda\leq\lambda_{1}$ and all $v\in\mathcal{K}\cap\partial
B_{\rho}(0)$,
\begin{align*}
Tv  &  \leq\phi^{-1}\left(  \int_{a}^{b}\lambda m(x)f(v)+\mu
n(x)g(v)dx\right)  \delta_{\Omega}\\
&  \leq\phi^{-1}\left(  \lambda C\int_{a}^{b}m(x)dx+c_{1}\mu\int_{a}%
^{b}n(x)v^{r_{1}}dx\right)  \delta_{\Omega}\\
&  \leq\phi^{-1}\left(  \lambda_{1}MC+\mu N\rho^{r_{1}}\right)  \delta
_{\Omega}\\
&  \leq\rho\text{ in }\Omega.
\end{align*}
This tells us that $\left\Vert Tv\right\Vert _{\infty}\leq\rho=\left\Vert
v\right\Vert _{\infty}$ for all $v\in\mathcal{K}\cap\partial B_{\rho}(0).$

Thus, Lemma \ref{teorema de punto fijo de conos en expansion y compresion}
says that $T$ has a fixed point in $\mathcal{K}\cap(\overline{B_{R}%
(0)}\setminus B_{\rho}(0))$. \qed

\subsection{Proof of item (III)}

\textit{Proof of Theorem 1.1 (III)}. In order to prove (III) we combine Lemma
\ref{lema subsol del problema concavo convexo} and the inequality
\eqref{desigualdad para probar inciso iii del concavo convexo}. Let
$0<\lambda<\Lambda.$ By the definition of $\Lambda$ there exists
$\overline{\lambda}\in\left(  \lambda,\Lambda\right]  $ and $u_{\overline
{\lambda}}\in\mathcal{P}^{\circ}$ solution of \eqref{problema concavo convexo}
associated to $\overline{\lambda}$. Since $\lambda<\overline{\lambda}$ it
follows that $u_{\overline{\lambda}}$ is a supersolution
\eqref{problema concavo convexo} associated to $\lambda$. Now, thanks to Lemma
\ref{lema subsol del problema concavo convexo} there exists $\varepsilon>0$
such that $v=\mathcal{S}_{\phi}(\varepsilon m\delta_{\Omega}^{q})$ is a
subsolution of \eqref{problema concavo convexo} associated to $\lambda$.
Moreover, taking $\varepsilon$ smaller if necessary, we get that $v\leq
u_{\overline{\lambda}}$. Now, (III) follows from Lemma 2.2.\qed

\section{Comments about the hypothesis}

\label{section comments}

Let us introduce some concepts about Orlicz spaces indices. Given a
nonbounded, increasing, continuous function $\phi:[0,\infty)\rightarrow
\lbrack0,\infty)$ with $\phi(0)=0$, we define
\[
M(t,\phi):=\sup_{x>0}\frac{\phi(tx)}{\phi(x)}.
\]
This function is nondecreasing and submultiplicative with $M(1,\phi)=1$. Then,
thanks to e.g. \cite[Chapter 11]{maligranda}, the following limits exist:
\[
\alpha_{\phi}:=\lim_{t\rightarrow0^{+}}\frac{\ln M(t,\phi)}{\ln t},\quad
\beta_{\phi}:=\lim_{t\rightarrow\infty}\frac{\ln M(t,\phi)}{\ln t},
\]
and moreover, $0\leq\alpha_{\phi}\leq\beta_{\phi}\leq\infty$. These numbers
are called \textbf{Orlicz space indices} or \textbf{Matuszewska-Orlicz's
indices}, who introduced them in \cite{ma}.

As usual, we say that $\phi$ satisfies the $\Delta_{2}$ condition if there
exists $k>0$ such that
\[
\phi(2x)\leq k\phi(x)\quad\text{for all }x\geq0.
\]

\begin{remark}
\label{remark about indeces}\strut

\begin{enumerate}
\item[(i)] For $\varepsilon>0$, there exists $t_{1}>0$ such that $\phi(tx)\leq
t^{\alpha_{\phi}-\varepsilon}\phi(x)$ for all $x>0$ and $t\in\lbrack0,t_{1}]$.

\item[(ii)] Suppose that $\beta_{\phi}<\infty$. Then, for $\varepsilon>0$,
there exists $t_{2}>0$ such that $\phi(tx)\leq t^{\beta_{\phi}+\varepsilon
}\phi(x)$ for all $x>0$ and $t\in\lbrack t_{2},\infty)$. So, if $\beta_{\phi
}<\infty$ then $\phi$ satisfies the $\Delta_{2}$ condition.

\item[(iii)] If $x^{-p}\phi(x)$ is nondecreasing for all $x>0$, then
$\alpha_{\phi}\geq p$.

\item[(iv)] If $x^{-p}\phi(x)$ is nonincreasing for all $x>0$, then
$\beta_{\phi}\leq p$.

\item[(v)] The following relationships between the Orlicz space indices of
$\phi$ and $\phi^{-1}$ hold:
\[
\beta_{\phi}=\frac{1}{\alpha_{\phi^{-1}}}\quad\text{and}\quad\alpha_{\phi
}=\frac{1}{\beta_{\phi^{-1}}}.
\]
As usual, we set $1/0=\infty$ and $1/\infty=0$.\smallskip
\end{enumerate}
\end{remark}

We shall need the next two useful lemmas to prove Theorem
\ref{equivalencias de phi y phi'} below.

\begin{lemma}
[\cite{gustavsson_peetre}, page 34]%
\label{lemma alpha positivo y beta finito implican cotas con homeos} If
$0<\alpha_{\phi}\leq\beta_{\phi}<\infty$ then there exist $C,p,q>0$ such that
\[
C^{-1}\min\{t^{p},t^{q}\}\phi(x)\leq\phi(tx)\leq C\max\{t^{p},t^{q}%
\}\phi(x)\quad\text{for all }t,x\geq0.
\]

\end{lemma}

\begin{lemma}
[\cite{maligranda}, Theorem 11.7]\label{lemma beta phi finito sii delta2} The
function $\phi$ satisfies the $\Delta_{2}$ condition if and only if the
constant $\beta_{\phi}$ is finite.
\end{lemma}

\begin{theorem}
\label{equivalencias de phi y phi'} The following hypothesis for $\phi$ are equivalent:

\begin{enumerate}
\item[(i)] $0<\alpha_{\phi}\le\beta_{\phi} <\infty. $

\item[(ii)] $(\Phi) $.

\item[(iii)] $(\Phi^{\prime}) $.
\end{enumerate}
\end{theorem}

\textit{Proof}. It is obvious that (ii) implies (iii), and Lemma
\ref{lemma alpha positivo y beta finito implican cotas con homeos} shows that
(i) implies (ii). Let us prove that (iii) implies (i).

Since $\alpha_{\phi}=1/\beta_{\phi^{-1}}$, Lemma
\ref{lemma beta phi finito sii delta2} and Remark \ref{remark about indeces}
(v) tell us that $\alpha_{\phi}>0$ if and only if $\phi^{-1}$ satisfies
$\Delta_{2}$. Let us check that the first inequality in $(\Phi^{\prime})$
implies that $\phi^{-1}$ satisfies $\Delta_{2}.$ Indeed, taking into account
that
\[
\psi_{1}(t)\phi(x)\leq\phi(xt)\quad\text{for all }t,x>0,
\]
setting $y=\phi(x)$ and $s=\psi(t)$ we get that
\[
sy\leq\phi(\psi_{1}^{-1}(s)\phi^{-1}(y))\quad\text{for all }s,y>0.
\]
Since $\phi^{-1}$ is increasing its follows that
\[
\phi^{-1}(sy)\leq\psi_{1}^{-1}(s)\phi^{-1}(y)\quad\text{for all }s,y>0.
\]
This implies that $\phi^{-1}$ satisfies $\Delta_{2}$. Thus, $\alpha_{\phi}>0$.
Moreover, the second inequality in $(\Phi^{\prime})$ implies that $\phi$
satisfies $\Delta_{2}$. Then, $\beta_{\phi}<\infty$. \qed\medskip

The following two lemmas will be useful to compare the indices $\alpha_{\phi}$
and $\beta_{\phi}$ with our hypotheses (F), (G1) and (G2) stated in Section 1.

\begin{lemma}
\label{teorema que relaciona limites con los indices de orlicz}Let $q>0$.

\begin{enumerate}
\item[(i)] If $\displaystyle\lim_{t\rightarrow0^{+}}\frac{t^{q}}{\phi(t)}=0$
then $\alpha_{\phi}\leq q$.

\item[(ii)] If $\displaystyle\lim_{t\rightarrow\infty}\frac{t^{q}}{\phi(t)}=0$
then $\beta_{\phi}\geq q$.

\item[(iii)] If $\displaystyle\lim_{t\rightarrow0^{+}}\frac{t^{q}}{\phi
(t)}=\infty$ then $\beta_{\phi}\geq q$.

\item[(iv)] If $\displaystyle\lim_{t\rightarrow\infty}\frac{t^{q}}{\phi
(t)}=\infty$ then $\alpha_{\phi}\leq q$.
\end{enumerate}
\end{lemma}

\textit{Proof}. We start proving (i). If $\alpha_{\phi}>q$, by Remark
\ref{remark about indeces} (i) there exists $t_{1}>0$ such that
\[
\phi(tx)\leq t^{q}\phi(x)\quad\text{for all }x>0\ \text{and}\ t\in(0,t_{1}).
\]
Let us set $C=\phi(1)^{-1}$ and fix $x=1$. Using the above inequality we have
that $C\leq\frac{t^{q}}{\phi(t)}$ for all $t\in(0,t_{1})$, which contradicts
that $\lim_{t\rightarrow0^{+}}\frac{t^{q}}{\phi(t)}=0.$ Therefore, we must
have $\alpha_{\phi}\leq q$. Item (ii) follows similarly. Indeed, if
$\beta_{\phi}<q$, by Remark \ref{remark about indeces} (ii) we have that there
exists $t_{1}>0$ such that
\[
\phi(tx)\leq t^{q}\phi(x)\quad\text{for all }x>0\ \text{and }t>t_{1}\text{.}%
\]
We now again define $C=\phi(1)^{-1}$ and fix $x=1$. Employing the above
inequality we have that $C\leq\frac{t^{q}}{\phi(t)}$ for all $t>t_{1}$,
contradicting that $\lim_{t\rightarrow\infty}\frac{t^{q}}{\phi(t)}=0.$ Thus,
$\beta_{\phi}\geq q.$

We prove (iii). We notice first that
\begin{equation}
\lim_{t\rightarrow0^{+}}\frac{t^{q}}{\phi(t)}=\infty\quad\text{if and only
if}\quad\lim_{t\rightarrow0^{+}}\frac{t^{1/q}}{\phi^{-1}(t)}=0.\label{u}%
\end{equation}
Indeed, the first limit is true if for every sequence $\{t_{k}\}$ with $0<$
$t_{k}\rightarrow0$, it holds that $\frac{t_{k}^{q}}{\phi(t_{k})}%
\rightarrow\infty$. Thus, taking $s_{k}=\phi(t_{k})$ we have that
$0<s_{k}\rightarrow0$ and $\frac{\left[  \phi^{-1}(s_{k})\right]  ^{q}}{s_{k}%
}\rightarrow\infty$. Since $h(t)=t^{1/q}$ is continuous and converges to
$\infty$ as $t\rightarrow\infty$, it follows that $\frac{\phi^{-1}(s_{k}%
)}{s_{k}^{1/q}}\rightarrow\infty$, which is equivalent to $\frac{s_{k}^{1/q}%
}{\phi^{-1}(s_{k})}\rightarrow0$. Since $0\leq\frac{t^{1/q}}{\phi^{-1}(t)}$
for all $t>0$ it follows that $\lim_{t\rightarrow0^{+}}\frac{t^{1/q}}%
{\phi^{-1}(t)}=0$. Now, from (\ref{u}) and item (i) we deduce that
$\alpha_{\phi^{-1}}\leq1/q$, and recalling Remark \ref{remark about indeces}
(v) we get that $\beta_{\phi}\geq q$, and (iii) holds. Analogously, (iv)
follows from (ii), taking into account that
\[
\displaystyle\lim_{t\rightarrow\infty}\frac{t^{q}}{\phi(t)}=\infty
\quad\text{if and only if}\quad\lim_{t\rightarrow\infty}\frac{t^{1/q}}%
{\phi^{-1}(t)}=0,
\]
and using again Remark \ref{remark about indeces} (v). \qed

\begin{lemma}
\label{teorema que relaciona limites con los indices de orlicz ii} Let
$\phi:[0,\infty)\rightarrow\lbrack0,\infty)$ be a nonbounded, increasing,
continuous function with $\phi(0)=0$.

\begin{enumerate}
\item[(i)] If $q<\alpha_{\phi}$ then $\displaystyle \lim_{t\rightarrow0^{+}}
\frac{t^{q}}{\phi(t)} = \infty. $

\item[(ii)] If $q>\beta_{\phi}$ then $\displaystyle \lim_{t\rightarrow\infty}
\frac{t^{q}}{\phi(t)} = \infty. $

\item[(iii)] If $q<\alpha_{\phi}$ then $\displaystyle \lim_{t\rightarrow
\infty} \frac{t^{q}}{\phi(t)} = 0. $

\item[(iv)] If $q>\beta_{\phi}$ then $\displaystyle \lim_{t\rightarrow0^{+}}
\frac{t^{q}}{\phi(t)} = 0. $
\end{enumerate}
\end{lemma}


Let us note that the reciprocals of items (i) and (ii) of the above lemma
\textit{are not true}, see Example (e.1) below.\smallskip

\textit{Proof.} Let us begin by proving (i). Let $\varepsilon>0$ such that
$\alpha_{\phi}-\varepsilon>q$. By Remark \ref{remark about indeces} (i) there
exists $t_{1}>0$ such that $\phi(tx)\leq t^{\alpha_{\phi}-\varepsilon}\phi(x)$
for all$\ x>0$ and $t<t_{1}$. Taking $x=1$ we get that $\frac{1}%
{t^{\alpha_{\phi}-\varepsilon}}\leq\frac{\phi(1)}{\phi(t)}$ for $t<t_{1}$.
Multiplying by $t^{q}$ on both sides and taking limit as $t\rightarrow0^{+}$
it follows that
\[
\lim_{t\rightarrow0^{+}}\frac{t^{q}}{t^{\alpha_{\phi}-\varepsilon}}\leq
\lim_{t\rightarrow0^{+}}\frac{\phi(1)t^{q}}{\phi(t)}.
\]
Since $q<\alpha_{\phi}-\varepsilon$, the first limit is infinite, and so also
the second one. Thus, (i) is proved.

Analogously, let $\varepsilon>0$ such that $\beta_{\phi}+\varepsilon<q$. By
Remark \ref{remark about indeces} (ii) there exists $t_{1}>0$ such that
$\phi(tx)\leq t^{\beta_{\phi}-\varepsilon}\phi(x)$ for all$\ x>0$ and
$t>t_{1}$. Taking $x=1$ we have $\frac{1}{t^{\beta_{\phi}-\varepsilon}}%
\leq\frac{\phi(1)}{\phi(t)}$ for $t<t_{1}$. Multiplying by $t^{q}$ on both
sides and taking limit as $t\rightarrow\infty$ we get
\[
\lim_{t\rightarrow\infty}\frac{t^{q}}{t^{\beta_{\phi}+\varepsilon}}\leq
\lim_{t\rightarrow\infty}\frac{\phi(1)t^{q}}{\phi(t)}.
\]
Since $q>\beta_{\phi}+\varepsilon$, the first limit is infinite, and thus also
the second one.

On other hand, (iii) follows from (ii) noting that
\[
\displaystyle\lim_{t\rightarrow\infty}\frac{t^{q}}{\phi(t)}=0\quad\text{if and
only if}\quad\lim_{t\rightarrow\infty}\frac{t^{1/q}}{\phi^{-1}(t)}=\infty,
\]
and taking into account that $\alpha_{\phi}>q$ if and only if $\beta
_{\phi^{-1}}<1/q$. Similarly, (iv) follows from (i) noting that
\[
\displaystyle\lim_{t\rightarrow0^{+}}\frac{t^{q}}{\phi(t)}=0\quad\text{if and
only if}\quad\lim_{t\rightarrow0^{+}}\frac{t^{1/q}}{\phi^{-1}(t)}=\infty,
\]
and recalling that $\beta_{\phi}<q$ if and only if $\alpha_{\phi^{-1}}>1/q$.
\qed

\begin{corollary}
Let $q$, $r_{1}$ and $r_{2}$ be given by (F), (G1) and (G2)
respectively.\strut

\begin{enumerate}
\item Suppose that $\alpha_{\phi} $ is positive.

\begin{enumerate}
\item If $q<\alpha_{\phi} $ then the limit in (F) holds.
\end{enumerate}

\item Suppose that $\beta_{\phi} $ is finite.

\begin{enumerate}
\item If $r_{1} > \beta_{\phi}$ then the limit in (G1) holds.

\item If $r_{2} > \beta_{\phi} $ then the limit in (G2) holds.
\end{enumerate}
\end{enumerate}
\end{corollary}

\subsection{Examples}

Let us conclude the article with some examples of functions $\phi$. We suppose
$x\geq0$ and we extend the function oddly.

\begin{enumerate}
\item[a.] Let
\[
\phi(x)=x^{p_{1}}+x^{p_{2}},\quad\text{with}\ p_{1}\geq p_{2}>0.
\]
Since $\phi(x)/x^{p_{1}}$ is nonincreasing and $\phi(x)/x^{p_{2}}$ is
nondecreasing, we see that $\beta_{\phi}<\infty$ and $\alpha_{\phi}>0$.

\item[b.] Let
\[
\phi(x)=\frac{x^{p_{1}}}{1+x^{p_{2}}},\quad\text{with}\ p_{1}>p_{2}>0.
\]
Since $\phi(x)/x^{p_{1}}$ is nonincreasing and $\phi(x)/x^{p_{1}-p_{2}}$ is
nondecreasing, we get that $\beta_{\phi}<\infty$ and $\alpha_{\phi}>0$.

\item[c.] Let
\[
\phi(x)=x\left(  \left\vert \ln x\right\vert +1\right)  .
\]
We have that $\phi(x)/x^{2}$ is nonincreasing. Then, $\beta_{\phi}<\infty.$
Furthermore, given $p\in(0,1)$ there exists $T>0$ such that
\[
\phi(tx)\leq t^{p}\phi(x)\quad\text{for}\ t\in\lbrack0,T]\text{ and all }%
x\geq0.
\]
This inequality implies that $\alpha_{\phi}\geq1.$

\item[d.] Let
\[
\phi(x):=x-\ln(x+1).
\]
As in the above example, $\phi(x)/x^{2}$ is nonincreasing and then
$\beta_{\phi}<\infty.$ Also, there exist $C,T>0$ such that
\[
\phi(tx)\leq\ Ct\phi(x)\quad\text{for }t\in\lbrack0,T]\text{ and all }x\geq0.
\]
The above inequality implies that $\alpha_{\phi}\geq1.$ Moreover, since
\[
\lim_{t\rightarrow\infty}\frac{t^{q}}{\phi(t)}=\infty\quad\text{for all }q>1,
\]
thanks to Lemma \ref{teorema que relaciona limites con los indices de orlicz}
(iv) we deduce that $\alpha_{\phi}=1$.

\item[e.] Let $h:(0,\infty)\rightarrow(1,\infty)$ be an increasing
differentiable function such that $\lim_{t\rightarrow0^{+}}h(t)=1$,
\begin{equation}
\lim_{t\rightarrow\infty}\frac{qt^{q-1}h(t)}{h^{\prime}(t)}=\infty
\quad\text{for all}\ q>0,\label{limite de h en infinito}%
\end{equation}
and there exists $p_{1}>0$ such that
\begin{equation}
\lim_{t\rightarrow0^{+}}\frac{qt^{q-1}h(t)}{h^{\prime}(t)}=%
\begin{cases}
0 & \text{if }q>p_{1},\\
\infty & \text{if }q<p_{1}.
\end{cases}
\label{limite de h en cero}%
\end{equation}
Define
\[
\phi(x):=(\ln(h(x))^{p},\quad\text{with}\ p>0.
\]
By \eqref{limite de h en infinito}, $\phi$ satisfies the limit in (G2).
Moreover, from Lemma
\ref{teorema que relaciona limites con los indices de orlicz} (iv) we can
deduce that $\alpha_{\phi}=0.$ Then $\phi$ \textit{does not satisfy} the
hypothesis $(\Phi)$ (and $(\Phi^{\prime})$) at the introduction. And since
\eqref{limite de h en cero} holds it follows that
\[
\lim_{t\rightarrow0^{+}}\frac{t^{q}}{\phi(t)}=%
\begin{cases}
0 & \text{if }q>pp_{1}.\\
\infty & \text{if }q<pp_{1}.
\end{cases}
\]
Therefore, $\phi$ satisfies the limits in (F) and (G1). Let us exhibit next a
few particular cases.

\begin{itemize}
\item[e.1] Let
\[
\phi(x):=(\ln(x+1))^{p},\quad\text{with}\ p>0.
\]
A few computations show that $h(x)=x+1$ satisfies
\eqref{limite de h en infinito} and \eqref{limite de h en cero}. Moreover, we
can see that $\phi(x)/x^{p}$ is nonincreasing and thus $\beta_{\phi}\leq p,$
and since
\[
\lim_{t\rightarrow0^{+}}\frac{t^{q}}{\ln(t+1)}=\infty\quad\text{for all }q<1,
\]
by Lemma \ref{teorema que relaciona limites con los indices de orlicz} it
follows that $\beta_{\phi}=p$. This shows that the reciprocals of the items
(i) and (ii) in Lemma
\ref{teorema que relaciona limites con los indices de orlicz ii} \textit{are
not true}.

\item[e.2] Let
\[
\phi(x):=\operatorname{arcsinh}(x)=\ln\left(  \sqrt{x^{2}+1}+x\right)  .
\]
One can see that $h(x)=\sqrt{x^{2}+1}+x$ satisfies
\eqref{limite de h en infinito} and \eqref{limite de h en cero}.

\item[e.3] Let
\[
\phi(x):=\ln(\ln(x+1)+1).
\]
One can verify that $h(x)=\ln(x+1)+1$ satisfies
\eqref{limite de h en infinito} and \eqref{limite de h en cero}.
\end{itemize}
\end{enumerate}

\end{document}